# A hybrid approach to solve the high-frequency Helmholtz equation with source singularity in smooth heterogeneous media

Jun Fang

*Department of Mathematics, University of California, Irvine*

Jianliang Qian

*Department of Mathematics, Michigan State University*

Leonardo Zepeda-Núñez

*Department of Mathematics, University of California, Irvine*

Hongkai Zhao

*Department of Mathematics, University of California, Irvine*

**Abstract**

We propose a hybrid approach to solve the high-frequency Helmholtz equation with point source terms in smooth heterogeneous media. The method is based on the ray-based finite element method (ray-FEM) [34], whose original version can not handle the singularity close to point sources accurately. This pitfall is addressed by combining the ray-FEM, which is used to compute the smooth far-field of the solution accurately, with a high-order asymptotic expansion close to the point source, which is used to properly capture the singularity of the solution in the near-field. The method requires a fixed number of grid points per wavelength to accurately represent the wave field with an asymptotic convergence rate of $\mathcal{O}(\omega^{-1/2})$, where $\omega$ is the frequency parameter in the Helmholtz equation. In addition, a fast sweeping-type preconditioner is used to solve the resulting linear system.

We present numerical examples in 2D to show both accuracy and efficiency of our method as the frequency increases. In particular, we provide numerical evidence of the convergence rate, and we show empirically that the overall complexity is $\mathcal{O}(\omega^2)$ up to a poly-logarithmic factor.

## 1. Introduction

The numerical solution of time-harmonic wave propagation in heterogeneous media is of paramount importance in a variety of applications such as medical imaging, oil exploration, nondestructive testing, noise reduction, radar and sonar technology.

In the constant density acoustic approximation the time harmonic wave propagation is modeled by the Helmholtz equation, which is given by

$$-\left(\Delta + \omega^2 m(\mathbf{x})\right) u(\mathbf{x}) = f(\mathbf{x}), \quad \mathbf{x} \in \Omega \subset \mathbb{R}^d \tag{1}$$

plus absorbing or radiation boundary conditions, where, $d$ is the dimension, $\omega$ is the angular frequency, $f$ is the source term, $m(\mathbf{x}) = 1/c^2(\mathbf{x})$ is the squared slowness, $c(\mathbf{x})$ is the wave speed, and $u$ is the unknown wave field to be computed, which becomes more oscillatory as the frequency increases.

Given the oscillatory behavior of $u$, we know, by the Nyquist-Shannon sampling criterion [64], that $\mathcal{O}(\omega^d)$ degrees of freedom are sufficient to accurately represent a wave field oscillating at frequency $\omega$. In addition, recent work [32] showed that this number of degrees of freedom is also intrinsic the Helmholtz equation, in the sense that the least dimension of a linear space that approximates the wavefields of a high-frequency Helmholtz equation to a given accuracy can grow with $\mathcal{O}(\omega^d)$ as $\omega \to \infty$. Hence, the optimal complexity to solve (1) is $\mathcal{O}(\omega^d)$, up to possible poly-log factors. However, computing the numerical solution of (1) in the high-frequency regime, i.e., when the number of degrees of freedom increases with $\omega$, accurately and efficiently is notoriously hard.

In view of the accuracy, the main challenge is rooted in the pollution error[1], which is ubiquitous in most piecewise discretizations. The main adverse effect of the pollution error manifests in large phase shifts in the propagating waves due to dispersion error of the discretization, even if the Nyquist-Shannon sampling criterion is respected, i.e., the mesh size, $h$, is comparable to the wavelength, $O(\omega^{-1})$. This implies that oversampling is required to obtain an accurate solution, leading to linear systems with suboptimal degrees of freedom [10, 9]. Several methods have been proposed to attenuate, and ultimately eliminate, the pollution error using either non-polynomial basis [56, 11, 12, 41], non-standard variational formulations [9, 69] or local polynomials refinements [8, 53, 55, 75, 83, 68]. Although some of the methods mentioned above are able to eliminate the pollution error, the resulting linear systems are either dense or extremely ill-conditioned thus non amenable to be solved fast, resulting in a suboptimal complexity [57, 39].

---
[1]The ratio between numerical error and best approximation error from a discrete finite element space is $\omega$ dependent.

In view of the computational efficiency, the main challenge is to solve the linear systems in quasi-linear time with respect to the number of unknowns. Standard sparse linear algebra algorithms based on nested dissection [37] and multi-frontal methods [26] have a suboptimal complexity, and they are prohibitively expensive memory-wise in dimension greater than two [25, 24, 3] even in the case of highly distributed codes [44], despite recent advances on compressed solvers [2, 73] using hierarchical compression techniques [13, 21]. Standard iterative methods[2] require a large, $\omega$-dependent, amount of iterations to converge, due to the indefinite character of the resulting linear systems, resulting in super-linear complexities with respect to the number of unknowns [33]. In this front many new highly efficient preconditioners have been proposed in the last 10 years that achieve the quasi-linear cost [29, 30, 65, 23, 72, 80]; however, most of them use low-order discretizations so that they require oversampling to produce accurate solutions, thus resulting in suboptimal complexities with respect to the frequency.

Even though the amount of literature dealing with both issues separately is vast [33, 79, 12, 41, 55, 83, 40, 66, 58, 31, 29, 65, 23, 72], only a few references deal with both issues *simultaneously*. We refer to, for example, [68] in which a hybridizable discontinuous Galerkin methods is coupled with the method of polarized traces, [81, 45] in which an integral version of the Helmholtz equation is coupled with sparsification and a fast preconditioner, and [34] in which an adaptive discretization is built by learning the dominant wave directions. The last approach is an important building block in the current paper and it will be explained in the sequel.

Besides the issues mentioned before, the accuracy and convergence rate of the discretization degrades greatly when the source in the right-hand side of (1) presents singularities [76, 51, 50]. Many applications model a point source using a Dirac delta for $f$, thus introducing a distributional right-hand side. In such case, mathematically, a wavefield has singularities in both the amplitude and the phase at the location of source. If the singularities are not properly resolved they will lead to significant errors propagating to the full domain, and, consequently, leading to the detriment of the accuracy of the solver. It is well known [5, 4, 19] that standard error estimates for finite elements are not valid when the source term is a distribution. Thus, some special treatments have been developed to handle such right-hand sides, mostly for elliptic problems, requiring specially tuned corrections (*cf.* [38, 28] and references therein), which are not easily applicable to the Helmholtz equation.

---

[2]Such as standard matrix splitting methods, including Jacobi, Gauss-Seidel and over-relaxed iterations (see [62] for further details); non-preconditioned Krylov space methods, including GMRES [63], Bi-CGSTAB [70], CG, etc; and standard multi-grid methods [35].

In the present paper we develop a hybrid approach to solve the high-frequency Helmholtz equation (1) with source singularity in smooth heterogeneous media efficiently. In particular, the method achieves the following:

- it only requires a fixed number of grid points per wavelength to accurately represent the solution with an overall linear computational complexity up to poly-log factors,

- it possess an asymptotic convergence rate of $\mathcal{O}(\omega^{-\frac{1}{2}})$, and

- it can handle sources modeled by Dirac delta seamlessly.

The hybrid approach presented in this paper is a natural extension of the ray-FEM [34], combining Babich's expansion [7, 61] to properly address the drawbacks of the former. The ray-FEM method is able to handle high-frequency problems accurately in quasi-linear time with respect to the intrinsic degrees of freedom and it has a convergence rate of $\mathcal{O}(\omega^{-\frac{1}{2}})$. However, as it will be explained in the sequel, it can not handle problems with point sources modeled by a Dirac delta. In such a scenario, we use the machinery developed in [61] to accurately remove the singularity from the right-hand side replacing it for a smooth one. The resulting equation is then solved using the ray-FEM method.

The properties of the hybrid approach presented above are achieved by complementing the two methods in the following manner:

- the Babich's expansion provides an efficient algorithm to obtain an accurate solution close to the singularity where the ray-FEM is inaccurate, and,

- the ray-FEM provides an efficient algorithm to compute an accurate solution far from the singularity, in which caustics and turning waves can occur and thus render difficulties for the Babich's expansion.

**Remark 1.** *We point out that in the high-frequency regime, the notion of convergence is different. We aim to study the asymptotic behavior of the error $\|u_\omega - u_\omega^h\|$ as $\omega$ grows, where $u_\omega$ and $u_\omega^h$ are the analytical and numerical solution to* (1) *respectively. In practice, most finite elements estimates for the error are bounded by an expression of the form $C(\omega)(\omega h)^p$, for some $p$ that depends on the order of the method. Usually, the pollution error manifests itself as a growing dependence of $C(\omega)$* with respect to $\omega$. Even though some methods have been able to provide a $C(\omega)$ which is constant with respect to $\omega$ by increasing the polynomial order of the approximation, the numerical error does not decay as $\omega \to \infty$ if $\omega h$ is fixed. In our case, however, the constant of the method decreases as $\mathcal{O}(\omega^{-1/2})$ as $\omega$ increases. To the best of our knowledge this is the first method to achieve, at least empirically, such scaling for a high-frequency problem with a distributional right-hand side.

1.1. Motivation

In [34] the authors presented a ray-based finite element method (ray-FEM) to solve (1) in smooth media based on a procedure to adaptively learn the basis functions to represent the wave field specific to the underlying medium and source. The approach was motivated by the geometric optics ansatz, in which the solution can take the form

$$u(\mathbf{x}) = \text{superposition of } \{A_n(\mathbf{x})e^{i\omega\phi_n(\mathbf{x})}\}_{n=1}^{N(\mathbf{x})} + O\left(\frac{1}{\omega}\right), \quad (2)$$

where $N(\mathbf{x})$ is the number of fronts/rays passing through $\mathbf{x}$, the phases $\phi_n$ are independent of the frequency $\omega$ and the amplitudes $A_n$ are mildly dependent on the frequency $\omega$. Furthermore, the phase functions can be linearly approximated in the form

$$\phi_n(\mathbf{x}) \approx \phi_n(\mathbf{x}_0) + \nabla\phi_n(\mathbf{x}_0) \cdot (\mathbf{x} - \mathbf{x}_0) = \phi_n(\mathbf{x}_0) + |\nabla\phi_n(\mathbf{x}_0)|\widehat{\mathbf{d}}_n(\mathbf{x}_0) \cdot (\mathbf{x} - \mathbf{x}_0), \quad (3)$$

where $\widehat{\mathbf{d}}_n(\mathbf{x}_0) := \frac{\nabla\phi_n(\mathbf{x})}{|\nabla\phi_n(\mathbf{x})|}$ are called the ray directions [17] or the dominant wave directions [20].

In [34], the local dominant plane-wave directions are learned/extracted first from a low-frequency wave, which probes the same medium using the same source. Those local dominant plane-wave directions are used to build local basis composed of polynomials[3] modulated by plane-waves. The resulting algorithm only requires a fixed number of grid points per wavelength to achieve both stability and convergence without oversampling. Moreover, a fast solver was developed for solving the resulting linear system with linear complexity up to poly-log factors. The method can efficiently compute wave fields when the source is far away, even in the presence of caustics. However, the ray-FEM cannot handle the singularities at point sources for two reasons:

---

[3]In this work, for simplicity we used a low degree polynomial, it is possible to use higher degree basis; however, the asymptotic behavior would remain the same.

- the traditional geometrical-optics amplitude at the source location is singular, which indicates that the geometric optics ansatz breaks down at the source points and it is difficult to handle for a finite element based method;

- the phase is also singular, i.e., the curvature of the circular wave-front goes to infinity at the source location, which makes the ray direction extraction, such as the numerical micro-local analysis (NMLA) [15, 16, 17], infeasible.

On the other hand, the Babich's expansion [7], which is a Hankel-based asymptotic expansion, can capture source singularity and overcome the above difficulties near the source in heterogeneous media. The components of the expansion can be numerically computed by high-order Eulerian asymptotic methods [61] to yield accurate solutions in the neighborhood of the point source.

The reasons for preferring the Babich's ansatz rather than the usual geometrical-optics ansatz are well illustrated in [61] and are briefly summarized below. If we apply the usual asymptotic expansion of the solution for the Helmholtz equation of a point source in an inhomogeneous medium, then we end up with the following systems:

$$u(\mathbf{x}, \mathbf{x}_0) = e^{i\omega\tau} \sum_{s=0}^{\infty} A_s(\mathbf{x}, \mathbf{x}_0) \frac{1}{(i\omega)^{s-\frac{(d-1)}{2}}}, \qquad (4)$$

where $\tau = \tau(\mathbf{x}, \mathbf{x}_0)$ is the phase satisfying the eikonal equation

$$\nabla\tau \cdot \nabla\tau = m^2(\mathbf{x}), \quad \tau(\mathbf{x}_0, \mathbf{x}_0) = 0, \qquad (5)$$

and $A_s = A_s(\mathbf{x}, \mathbf{x}_0)$ satisfy a recursive system of transport equations along rays,

$$2\nabla\tau \cdot \nabla A_s + A_s \Delta\tau = -\Delta A_{s-1}, \quad s = 0, 1, \cdots, \quad A_{-1} \equiv 0. \qquad (6)$$

(Note: I changed **r** by **x**, to be more consistent with the notation above.)

Alas, a difficulty arises immediately: how to initialize $A_s$ at the source point for this system of equations. In addition, when $d$ is even, the ray series (4) does not yield a uniform asymptotic form close to the source. When $d = 3$, Avila and Keller [6] were able to find the initial data for $A_s$ using the boundary layer method, but the case of $d = 2$ was left incomplete. In practice, such difficulties in initializing amplitudes are handled in an *ad hoc* manner. The amplitudes are initialized a little bit away from the point source using amplitudes for a medium with constant refractive index [71, 60, 47, 49, 48]; consequently, the resulting numerical asymptotic solution is not uniform near the source.

To overcome these initialization difficulties Babich [7] proposed an asymptotic series defined by the first Hankel function as an ansatz to expand the underlying highly-oscillatory wave field; the equation for the phase is the usual eikonal equation, but the resulting transport equations for the amplitudes are easily initialized. Moreover, Babich's expansion ensures that the Hankel-based ansatz yields a uniform asymptotic solution as $\omega \to \infty$ within the neighborhood of the point source and away from it. However, the Babich's ansatz based on the Eikonal equation for the phase function and transport equations for amplitudes may encounter difficulties due to caustics and turning/crossing rays.

The hybrid approach developed in this paper is aimed to combine the strengths of the ray-FEM and Babich's expansion to solve the high-frequency Helmholtz equation with source singularity. Numerical examples in 2D demonstrate that the proposed hybrid method achieves an asymptotic convergence rate of $\mathcal{O}(\omega^{-\frac{1}{2}})$ with fixed number of grid points per wavelength and a total empirical complexity of $\mathcal{O}(\omega^2)$ up to a poly-logarithmic factor.

1.2. *The hybrid approach*

We decompose $u$, the solution to (1) with point source term $f(\mathbf{x}) = \delta(\mathbf{x} - \mathbf{x}_0)$, into two components:
$$u(\mathbf{x}) = u_{\text{near}}(\mathbf{x}) + u_{\text{far}}(\mathbf{x}), \tag{7}$$
where $u_{\text{near}}$ is the near-field solution, which captures the source singularity, and $u_{\text{far}}$ is the far-field solution, which is highly-oscillatory but does not contain singularities. We insert (7) into (1) and we have that
$$-\left(\Delta + \omega^2 m\right) u_{\text{far}} = \delta(\mathbf{x} - \mathbf{x}_0) + \left(\Delta + \omega^2 m\right) u_{\text{near}}. \tag{8}$$
Moreover, we suppose that $u_{\text{near}}$ has the form
$$u_{\text{near}}(\mathbf{x}) = u_b(\mathbf{x})\chi_\epsilon(\mathbf{x}), \tag{9}$$
where $u_b$ is the approximation given by the Babich's expansion and $\chi_\epsilon$ is a smooth cut-off function satisfying,
$$\chi_\epsilon(\mathbf{x}) = \begin{cases} 1, & \text{if } |\mathbf{x} - \mathbf{x}_0| < \epsilon, \\ 0, & \text{if } |\mathbf{x} - \mathbf{x}_0| > 2\epsilon, \end{cases} \tag{10}$$
where $\epsilon$ is a fixed small number such that there are no caustics or ray crossing in $|\mathbf{x} - \mathbf{x}_0| < 2\epsilon$. Following a standard computation we have that
$$\left(\Delta + \omega^2 m\right) u_{\text{near}} = \left(\Delta u_b + \omega^2 m u_b\right)\chi_\epsilon + 2\nabla u_b \cdot \nabla \chi_\epsilon + u_b \Delta \chi_\epsilon, \tag{11}$$

furthermore, given that $\Delta u_b + \omega^2 m u_b = -\delta(\mathbf{x} - \mathbf{x}_0)$ and plugging into (8), we have that
$$-\left(\Delta + \omega^2 m\right) u_{\text{far}} = 2\nabla u_b \cdot \nabla \chi_\epsilon + u_b \Delta \chi_\epsilon, \tag{12}$$
whose right-hand side is smooth. Moreover, it can be easily computed accurately: $\chi_\epsilon$ and its derivatives are known analytically, and it is possible to compute $u_b$ accurately and efficiently in the support of $\chi_\epsilon$ using the method developed in [61], as will be reviewed briefly in the sequel.

Given that $u_{\text{far}} = u - u_{\text{near}} \simeq (1 - \chi_\epsilon)u$, i.e., $u_{\text{far}}$ is the far-field solution of the Helmholtz equation, it satisfies absorbing or radiation conditions. Then we can solve the equation (12) using the techniques developed in [34] for high-frequency Helmholtz equation in smooth media, which has a total complexity proportional to $\omega^2$ up to polylogarithmic factors.

The proposed algorithm to compute the solution to (1) with a point source can be distilled to the following steps:

- compute the asymptotic solution $u_b$ given by the Babich's expansion in a neighborhood of the source point,

- build the right-hand side in (12),

- solve (12), using the ray-FEM method with the adaptive learning basis approach proposed in [34],

- add the near field part $u_{\text{near}}$ and the far field part $u_{\text{far}}$.

We provide details of each step in the sequel.

*1.3. Outline*

The outline of the paper is as the following: we give a brief introduction of the Babich's expansion in Section 2 and we review the adaptive learning ray-FEM method in Section 3. In Section 4, we provide the error analysis of the hybrid method in each step. In Section 5 we provide a description of fast linear solvers that we are using. The full algorithm is described in Section 6 and numerical results are presented in Section 7. Finally, conclusions are summarized in Section 8.

## 2. Babich's expansion

The reduction of (7) to (12) relies on computing a local approximation of the solution close to the source. This is achieved by using a specific asymptotic expansion of the solution, usually referred to as the Babich's expansion [7], which we briefly review here.

To solve (1) asymptotically when $\omega \to \infty$, Babich proposed the following Hankel-based ansatz [7] to expand the solution in a neighborhood of the source,

$$u(\mathbf{x}) = u_b(\mathbf{x}, \mathbf{x}_0, \omega) := \sum_{p=0}^{\infty} v_p(\mathbf{x}, \mathbf{x}_0) f_{p-(d-2)/2}(\omega, \phi(\mathbf{x}, \mathbf{x}_0)), \tag{13}$$

where $d$ is the dimension,

$$f_q(\omega, \xi) := i\frac{\sqrt{\pi}}{2} e^{iq\pi} \left(\frac{2\xi}{\omega}\right)^q H_q^{(1)}(\omega\xi), \tag{14}$$

and $\phi$ is the phase or travel time function satisfying the eikonal equation,

$$|\nabla \phi(\mathbf{x}, \mathbf{x}_0)| = \frac{1}{c(\mathbf{x})}, \tag{15}$$

and $\{v_p\}_{p=0,1,2,\ldots}$ are assumed to be smooth functions satisfying a recursive system of transport equations,

$$\nabla \phi^2(\mathbf{x}) \cdot \nabla v_p(\mathbf{x}) + \left[(2p-d)m(\mathbf{x}) + \frac{1}{2}\Delta\phi^2(\mathbf{x})\right] v_p(\mathbf{x}) = \frac{1}{2}\Delta v_{p-1}(\mathbf{x}) \tag{16}$$

with initial conditions

$$v_{-1} \equiv 0, \quad v_0(\mathbf{x}_0) = \frac{m_0^{(d-2)/2}}{2\pi^{(d-1)/2}}, \quad m_0 = m(\mathbf{x}_0). \tag{17}$$

The problem is then reduced to computing each term of the expansion (13), defined as $u_{b_p} = v_p f_{p-(d-2)/2}$, by solving the the PDE systems (15) and (16). By using the high-order Eulerian asymptotic method developed in [61], it is possible obtain accurate numerical solutions even close to the source.

In the present work we only need to compute the first two terms of the Babich's expansion in order to obtain the desired convergence rate of $\mathcal{O}(\omega^{-\frac{1}{2}})$.

The algorithm to compute the expansion can be summarized as follows:

- firstly, the phase $\phi$ is computed using a fifth-order Lax-Friedrichs weighted non-

oscillatory (LxF-WENO) scheme [43, 82, 77] with a sixth-order factorization[4] [36, 47, 48] around the source;

- secondly, the first amplitude coefficient $v_0$ is computed using the third-order LxF-WENO scheme with the third-order factorization around the source;

- then, the second amplitude coefficient $v_1$ is computed by the third-order LxF-WENO scheme with the first-order factorization around the source;

- finally, the solution is approximated by replacing the numerically computed phase and amplitude coefficients into $u_{b_0} + u_{b_1}$.

Note that the different orders of the LxF-WENO schemes above are calibrated to each other for the different subproblems, the reason is well explained in Section 3 of [61]. Below is a summary of the algorithm to approximate the symptotic Babich's expansion $u_b$ and its gradient $\nabla u_b$ in the disk $D_{2\epsilon} := \{\mathbf{x} \in \Omega : |\mathbf{x} - \mathbf{x}_0| \leq 2\epsilon\}$ with mesh size $h$, and we denote by $u_b^h, \nabla u_b^h, u_{b_0}^h, u_{b_1}^h, \phi_h, v_0^h, v_1^h$ the numerically computed quantities of $u_b, \nabla u_b, u_{b_0}, u_{b_1}, \phi, v_0, v_1$ respectively. For further details, we refer the reader to [61].

---

**Algorithm 1** Babich's Expansion

1: **function** $[u_b^h, \nabla u_b^h, \{\boldsymbol{d}_b^h\}_{D_{2\epsilon}}] = Babich(\mathbf{x}_0, c, \omega, h)$
2:      $[\phi_h, \nabla\phi_h, v_0^h, \nabla v_0^h, v_1^h, \nabla v_1^h] = Eikonal\text{-}Transport(\mathbf{x}_0, c, h)$
3:      $\boldsymbol{d}_b^h = \frac{\nabla\phi_n(\mathbf{x})}{|\nabla\phi_n(\mathbf{x})|}$      ▷ Computing ray directions
4:      $coef_0 = i\frac{\sqrt{\pi}}{2}$, $coef_1 = i\frac{\sqrt{\pi}}{\omega}\exp(i\pi)$      ▷ Coefficients of $f_q$ in (14)
5:      $\tilde{f}_0^h = H_0^{(1)}(\omega\phi_h)$, $\tilde{f}_1^h = \phi_h H_1^{(1)}(\omega\phi_h)$      ▷ Hankel based terms of $f_q$ in (14)
6:      $u_{b_0}^h = coef_0 v_0^h \tilde{f}_0^h$, $u_{b_1}^h = coef_1 v_1^h \tilde{f}_1^h$      ▷ The first two terms in (13)
7:      $u_b^h = u_{b_0}^h + u_{b_1}^h$      ▷ Approximating the Babich's expansion
8:      $\nabla\tilde{f}_0^h = -\omega H_1^{(1)}(\omega\phi_h)\nabla\phi_h$      ▷ Computing gradients
9:      $\nabla\tilde{f}_1^h = 2H_1^{(1)}(\omega\phi_h)\nabla\phi_h - \omega\phi_h H_2^{(1)}(\omega\phi_h)\nabla\phi_h$
10:     $\nabla u_{b_0}^h = coef_0(\nabla v_0^h \tilde{f}_0^h + v_0^h \nabla\tilde{f}_0^h)$, $\nabla u_{b_1}^h = coef_1(\nabla v_1^h \tilde{f}_1^h + v_1^h \nabla\tilde{f}_1^h)$
11:     $\nabla u_b^h = \nabla u_{b_0}^h + \nabla u_{b_1}^h$      ▷ Approximating the gradient of Babich's expansion
12: **end function**

---

[4]The solution is represented as a product or sum of the analytical solution to a homogeneous medium and an unknown factor or perturbation, which is smooth.

**Algorithm 2** Eikonal/Transport Solver
―――――――――――――――――――――――――――――――――――――――――――――
1: **function** $[\phi_h, \nabla\phi_h, v_0^h, \nabla v_0^h, v_1^h, \nabla v_1^h] = \textit{Eikonal-Transport}(\mathbf{x}_0, c, h)$
2:     $\phi_h = \textit{LxF-WENO-Fac}(5, 5, \mathbf{x}_0, c, h)$         ▷ Computing the phase in (15)
3:     $\nabla\phi_h^2 = \textit{WENO}(3, \phi_h^2, h)$, $\Delta\phi_h^2 = \textit{FD}(4, \phi_h^2, h)$, $\nabla\phi_h = \nabla\phi_h^2/2\phi_h$
4:     $v_0^h = \textit{LxF-WENO-Fac}(3, 3, \nabla\phi_h^2, \Delta\phi_h^2, \mathbf{x}_0, c, h)$         ▷ Computing the first amplitude coefficient in (16)
5:     $\nabla v_0^h = \textit{FD}(4, v_0^h, h)$, $\Delta v_0^h = \textit{FD}(4, v_0^h, h)$
6:     $v_1^h = \textit{LxF-WENO-Fac}(3, 1, \nabla\phi_h^2, \Delta\phi_h^2, \Delta v_0^h, \mathbf{x}_0, c, h)$         ▷ Computing the second amplitude coefficient in (16)
7:     $\nabla v_1^h = \textit{FD}(4, v_1^h, h)$
8: **end function**
―――――――――――――――――――――――――――――――――――――――――――――

**Remark 2.** *In Algorithm 1, we compute $u_b^h$ and $\nabla u_b^h$ only in the annulus $\epsilon \leq |\mathbf{x} - \mathbf{x}_0| \leq 2\epsilon$ to build the right-hand side of equation (12) and $\{\boldsymbol{d}_b^h\}_{D_{2\epsilon}}$ to be used in the ray-FEM. In Algorithm 2, we denote by LxF-WENO-Fac$(p, q, \cdot)$ the p-th order LxF-WENO scheme with q-th order factorization around the source, WENO$(p, \cdot)$ the p-th order WENO scheme and FD$(p, \cdot)$ the p-th order finite difference scheme. $\nabla\phi_h^2$, $\nabla\phi_h$ and $\Delta\phi_h^2$ are computed except for the points around the source.*

## 3. Ray-FEM

As shown in [61], the Babich's expansion provides an accurate solution in the neighborhood of the source, provided that high-order methods are used to compute the phase and amplitude terms of the expansion. However, this approach cannot directly handle caustics or crossing rays away from the source. This is a consequence of computing the viscosity solution of the eikonal equation, which does not take into account possible multivalued solutions of the eikonal equation. Therefore, we only use this expansion in a region around the source without caustics occurring. Such a region always exists for the eikonal equation under consideration as shown in [6, 67]. In the far field, however, we compute the waves by using the ray-FEM method developed in [34], which is a stable finite element method without oversampling, specially adapted for wave propagation. The ray-FEM method is based on a geometric optic ansatz and its local approximation via a superposition of plane-waves propagating in a set of dominant directions. We use this method to solve equation (12).

### 3.1. Geometric optics ansatz

The geometric optics ansatz approximates the solution to the Helmholtz problem asymptotically in $\omega$, by a superposition of several wave-fronts of the form in (2). Moreover, except for a small set of points, e.g., source/focus points, caustics, and discontinuities of the medium, phases $\phi_n$ and amplitudes $A_n$ are single-valued functions satisfying the following PDE system,

$$\text{(eikonal)} \quad |\nabla \phi_n| = \frac{1}{c}, \qquad \text{(transport)} \quad 2\nabla \phi_n \cdot \nabla A_n + A_n \Delta \phi_n = 0, \qquad (18)$$

each defined in a suitable domain with suitable boundary conditions [14].

Based on the above ansatz, a local plane-wave approximation can be derived at any point where $\phi_n$ and $A_n$ are smooth with variations on an $\mathcal{O}(1)$ scale. Indeed, we define

$$\widehat{\mathbf{d}}_n := \frac{\nabla \phi_n(\mathbf{x}_0)}{|\nabla \phi_n(\mathbf{x}_0)|} = c(\mathbf{x}_0)\nabla \phi_n(\mathbf{x}_0) \qquad (19)$$

as the ray directions of the wave-fronts at an observation point $\mathbf{x}_0$, $k(\mathbf{x}) = \omega/c(\mathbf{x})$, and

$$B_n(\mathbf{x}) = (A_n(\mathbf{x}_0) + \nabla A_n(\mathbf{x}_0)(\mathbf{x} - \mathbf{x}_0))e^{i\omega(\phi_n(\mathbf{x}_0) - \nabla \phi(\mathbf{x}_0)\cdot\mathbf{x}_0)} \qquad (20)$$

the affine complex amplitude on a small neighborhood $|\mathbf{x} - \mathbf{x}_0| < h \ll 1$. By using a Taylor expansion within this neighborhood, we obtain a local plane-wave approximation of the following form for the $n$-th wave-front,

$$u_n(\mathbf{x}) := B_n(\mathbf{x})e^{ik(\mathbf{x}_0)\widehat{\mathbf{d}}_n \cdot \mathbf{x}} + \mathcal{O}\left(h^2 + \omega h^2 + \frac{1}{\omega}\right). \qquad (21)$$

This motivates us to construct local finite element basis with mesh size $h = \mathcal{O}(\omega^{-1})$, in which an affine function is multiplied by plane-waves oscillating along the dominant ray directions, resulting in a local approximation similar to (21) but with asymptotic error $\mathcal{O}(\omega^{-1})$.

### 3.2. Ray-based FEM formulation

We use finite element methods to compute the solution to (12) with PML boundary conditions [18]. For simplicity, we consider the rectangular domain $\Omega = (-L_x, L_x) \times$

$(-L_y, L_y)$ in 2D. We introduce

$$\delta_x(x) = \begin{cases} \frac{C}{\delta_{pml}} \left( \frac{x+L_x-\delta_{pml}}{\delta_{pml}} \right)^2, & \text{if } x \in (-L_x, -L_x + \delta_{pml}), \\ 0, & \text{if } x \in (-L_x + \delta_{pml}, L_x - \delta_{pml}), \\ \frac{C}{\delta_{pml}} \left( \frac{x-L_x+\delta_{pml}}{\delta_{pml}} \right)^2, & \text{if } x \in (L_x - \delta_{pml}, L_x), \end{cases} \quad (22)$$

and similarly for $\delta_y(y)$. Here $\delta_{pml}$ is typically around a couple of wavelengths, and $C$ is an appropriate positive absorption constant independent of $\omega$.

Then the equation (1) can be re-written [30] as

$$-\nabla \cdot (\mathbf{D}\nabla u) - \omega^2 m s_x s_y u = s_x s_y f \text{ in } \Omega, \quad u = 0 \text{ on } \partial\Omega \quad (23)$$

where $\mathbf{D} = \begin{bmatrix} s_y/s_x & 0 \\ 0 & s_x/s_y \end{bmatrix}$, $s_x = 1+i\sigma_x(x)/\omega$ and $s_y = 1+i\sigma_y(y)/\omega$ with quadratic coefficient functions $\sigma_x(x)$ and $\sigma_y(y)$. The standard weak formulation is given by

$$\text{Find } u \in H_0^1(\Omega), \text{ such that } \mathcal{B}(u,v) = \mathcal{F}(v), \quad \forall v \in H_0^1(\Omega), \quad (24)$$

where

$$\mathcal{B}(u,v) := \int_\Omega (\mathbf{D}\nabla u) \cdot \nabla \bar{v} dV - \omega^2 \int_\Omega m s_x s_y u \bar{v} dV \quad (25)$$

$$\mathcal{F}(v) := \int_\Omega s_x s_y f \bar{v} dV, \quad f = 2\nabla u_b \cdot \nabla \chi_\epsilon + u_b \Delta \chi_\epsilon. \quad (26)$$

The domain, $\Omega$, is discretized with a standard regular triangulated mesh, with mesh size $h$. The resulting mesh is denoted by $\mathcal{T}_h = \{K\}$, where $K$ represents a triangle element of the mesh.

We define the standard local approximation space by

$$V_S(K) = \text{span}\{\varphi_j(\mathbf{x}), \quad j = 1, 2, 3, \quad \forall \mathbf{x} \in K\}, \quad (27)$$

where $\{\varphi_j(\mathbf{x})\}_{j=1}^3$ are standard $\mathbb{P}1$ finite element nodal basis functions over the element $K$.

We modify the $\mathbb{P}1$ finite elements by incorporating the ray information and we define the ray-based local approximation space by

$$V_{Ray}(K) = \text{span}\{\varphi_j(\mathbf{x})e^{ik_j \widehat{\mathbf{d}}_{j,l} \cdot \mathbf{x}}, \quad j = 1, 2, 3, \quad l = 1, ..., n_j, \quad \forall \mathbf{x} \in K\}, \quad (28)$$

where we let $k_j = \omega/c(\mathbf{x}_j)$ be wave-numbers and $\{\widehat{\mathbf{d}}_{j,l}\}_{l=1}^{n_j}$ be $n_j$ ray directions at the vertex $V_j$ with coordinates $\mathbf{x}_j$.

Therefore, we can build the corresponding global approximation space and further the standard/ray-based FEM by replacing $H_0^1(\Omega)$ with $V_{S/Ray}(\mathcal{T}_h)$ in (24), where

$$V_{S/Ray}(\mathcal{T}_h) = \{v \in C^0(\overline{\Omega}) : v|_K \in V_{S/Ray}(K), \quad \forall K \in \mathcal{T}_h\}. \tag{29}$$

### 3.3. Adaptive learning of the basis

In the prequel we used the geometric optics ansatz to build an adaptive approximation space that incorporates ray information specific to the underlying Helmholtz equation. However, the ray directions, which depend on the medium and source distribution, are unknown quantities themselves, hence they need to be computed or estimated numerically.

Following [34] the dominant ray directions are learned by probing the same medium with the same source but using a relative low-frequency wave. That is, we solve the low-frequency Helmholtz equation

$$-\Delta \widetilde{u}(\mathbf{x}) - \widetilde{\omega}^2 m(\mathbf{x})\widetilde{u}(\mathbf{x}) = f(\mathbf{x}), \quad \mathbf{x} \in \Omega \subseteq \mathbb{R}^d, \tag{30}$$

plus suitable boundary (or radiation) conditions with the same $m(\mathbf{x}), f(\mathbf{x})$ but at a relative low-frequency $\widetilde{\omega} \sim \sqrt{\omega}$ on a mesh with size $h = \mathcal{O}(\widetilde{\omega}^{-2}) = \mathcal{O}(\omega^{-1})$ by using the standard FEM, which is quasi-optimal in that regime [54].

We then process the computed low-frequency wave field, with numerical microlocal analysis (NMLA), to extract the dominant ray directions locally in the far field. A brief introduction to NMLA can be found in Appendix A and [16, 17, 34]. However, NMLA cannot capture ray directions accurately near the source since the wave-fronts near the source have high curvatures. Instead we apply Algorithm 1 to obtain the ray directions near the source.

Once the local ray directions are extracted, we incorporate them into the ray-FEM space to solve the high-frequency Helmholtz equation (1). If necessary, ray estimation can be improved by iteratively applying NMLA to high-frequency wave field and then the resulting ray information can be incorporated into the ray-FEM space to obtain more accurate high-frequency waves.

## 4. Error analysis

In this section, we provide an asymptotic error estimate of the high-frequency solution $u$ to (1) by decomposing it into two parts: near field $u_{\text{near}}$ computed using Babich's expansion [61] and far field $u_{\text{far}}$ computed using the ray-FEM [34]. We remind the reader that in this section, we assume the wave speed $c$, phase function $\phi$, amplitude coefficients $v_1$, $v_2$ and amplitude function $A$ are smooth functions, and we only focus on the rectangular 2D domain with mesh size that scales as $\omega h = \mathcal{O}(1)$ in the high-frequency regime. We define notations $D_r := \{\mathbf{x} \in \Omega : |\mathbf{x} - \mathbf{x}_0| \leq r\}$, $D_{r_2 - r_1} := D_{r_2} \backslash D_{r_1} = \{\mathbf{x} \in \Omega : r_1 < |\mathbf{x} - \mathbf{x}_0| \leq r_2\}$.

## 4.1. Near-field solution: Babich's expansion

We recall the basic properties of the first kind Hankel functions [42, 74]

$$\frac{d}{dz} H_q^1(z) = q \frac{H_q^{(1)}(z)}{z} - H_{q+1}^{(1)}(z), \tag{31}$$

$$H_q^{(1)}(z) = \begin{cases} \mathcal{O}(z^{-1/2}), & \text{if } |z| \to \infty, \\ \mathcal{O}(\ln z), & \text{if } |z| \to 0 \text{ and } q = 0, \\ \mathcal{O}(z^{-q}), & \text{if } |z| \to 0 \text{ and } q \geq 1. \end{cases} \tag{32}$$

In the disk $D_{2\epsilon}$, $\phi(\mathbf{x}) \leq \mathcal{O}(\epsilon) \leq C_1 = $ constant, the Hankel based terms $f_q(\omega, \xi)$ have the following asymptotic form [7] for large $\omega$,

$$f_q(\omega, \xi) = \begin{cases} \mathcal{O}\left(\left(\frac{\xi}{\omega}\right)^q (\omega\xi)^{-1/2}\right) = \mathcal{O}(\omega^{-q-1/2} \xi^{q-1/2}), & \text{if } \omega\xi \geq C_2 = \text{constant}, \\ \mathcal{O}\left(\ln(\omega\xi)\right), & \text{if } \omega\xi \leq C_2 \text{ and } q = 0, \\ \mathcal{O}\left(\left(\frac{\xi}{\omega}\right)^q (\omega\xi)^{-q}\right) = \mathcal{O}(\omega^{-2q}), & \text{if } \omega\xi \leq C_2 \text{ and } q \geq 1. \end{cases} \tag{33}$$

Since the Babich's expansion (13) is approximated by the first two terms in Algorithm 1, i.e. $u_b \approx u_{b_0} + u_{b_1}$, where $u_{b_q} = v_q f_q$, $q = 0, 1$. We have the truncation error, asymptotically in $\omega$,

$$\|u_b - (u_{b_0} + u_{b_1})\|_{L^\infty(D_{2\epsilon})} = \mathcal{O}(\omega^{-5/2}). \tag{34}$$

On the other hand, the phase $\phi$ and amplitude coefficients $v_0, v_1$ are numerically computed by Algorithm 2. According to Theorem 5.1 in [61] and Remark 3 in [46], the $p$-th order LxF-WENO scheme combines with $q$-th order factorization for equations (15) and (16) yield $min(p,q)$-th order accuracy for smooth $\phi$ and $v$'s. Thus, we have

$$\|\phi - \phi_h\|_{L^\infty(D_{2\epsilon})} = \mathcal{O}(h^5), \quad \|v_0 - v_0^h\|_{L^\infty(D_{2\epsilon})} = \mathcal{O}(h^3), \quad \|v_1 - v_1^h\|_{L^\infty(D_{2\epsilon})} = \mathcal{O}(h). \tag{35}$$

Excluding a small neighborhood of the singular source point, i.e. $\mathbf{x} \in D_{2\epsilon - \eta}$, $\omega\phi(\mathbf{x}) \geq \mathcal{O}(\omega\eta) \geq C_2$ for a small positve number $\eta < \epsilon$, by the mean value theorem there exists $\varphi \geq \min\{\phi(\mathbf{x}), \phi^h(\mathbf{x})\} \geq \mathcal{O}(\eta)$ such that

$$|H_0^{(1)}(\omega\phi(\mathbf{x})) - H_0^{(1)}(\omega\phi^h(\mathbf{x}))| \begin{aligned} &= |-H_1^{(1)}(\omega\varphi)(\omega\phi(\mathbf{x}) - \omega\phi^h(\mathbf{x}))| \\ &\lesssim \mathcal{O}((\omega\varphi)^{-1/2}\omega h^5) \lesssim \mathcal{O}(\omega^{1/2}\eta^{-1/2}h^5), \end{aligned} \qquad (36)$$

where the constants for $\lesssim$ and $\mathcal{O}(\cdot)$ only depend on constants $C_1$ and $C_2$. Thus, we have

$$\|f_0 - f_0^h\|_{L^\infty(D_{2\epsilon-\eta})} \lesssim \|H_0^{(1)}(\omega\phi) - H_0^{(1)}(\omega\phi^h)\|_{L^\infty(D_{2\epsilon-\eta})} = \mathcal{O}(\omega^{1/2}\eta^{-1/2}h^5), \quad (37)$$

similarly, $\|f_1 - f_1^h\|_{L^\infty(D_{2\epsilon-\eta})} = \mathcal{O}(\omega^{-1/2}h^5 + \omega^{-3/2}\eta^{-1/2}h^5)$.

Hence,

$$\begin{aligned}
\|u_{b_0} - u_{b_0}^h\|_{L^\infty(D_{2\epsilon-\eta})} &= \|v_0 f_0 - v_0^h f_0^h\|_{L^\infty(D_{2\epsilon-\eta})} \\
&= \|(v_0 f_0 - v_0^h f_0) + (v_0^h f_0 - v_0^h f_0^h)\|_{L^\infty(D_{2\epsilon-\eta})} \\
&\leq \|v_0 - v_0^h\|_{L^\infty(D_{2\epsilon-\eta})}\|f_0\|_{L^\infty(D_{2\epsilon-\eta})} \\
&\quad + \|v_0^h\|_{L^\infty(D_{2\epsilon-\eta})}\|f_0 - f_0^h\|_{L^\infty(D_{2\epsilon-\eta})} \\
&= \mathcal{O}(\omega^{-1/2}\eta^{-1/2}h^3 + \omega^{1/2}\eta^{-1/2}h^5)
\end{aligned} \qquad (38)$$

and analogously,

$$\|u_{b_1} - u_{b_1}^h\|_{L^\infty(D_{2\epsilon-\eta})} = \mathcal{O}(\omega^{-3/2}h + \omega^{-1/2}h^5 + \omega^{-3/2}\eta^{-1/2}h^5) \qquad (39)$$

Therefore, under the assumption of $h = \mathcal{O}(\omega^{-1})$ and $\frac{C_2}{\omega} \leq \mathcal{O}(\eta) \leq \phi \leq \mathcal{O}(\epsilon) \leq C_1$, the asymptotic error with respect to $\omega$ of the Babich's expansion is

$$\begin{aligned}
\|u_b - u_b^h\|_{L^\infty(D_{2\epsilon-\eta})} &= \|u_b - u_{b_0}^h - u_{b_1}^h\|_{L^\infty(D_{2\epsilon-\eta})} \\
&= \|(u_b - u_{b_0} - u_{b_1}) + (u_{b_0} - u_{b_0}^h) + (u_{b_1} - u_{b_1}^h)\|_{L^\infty(D_{2\epsilon-\eta})} \\
&\leq \|u_b - u_{b_0} - u_{b_1}\|_{L^\infty(D_{2\epsilon-\eta})} \\
&\quad + \|u_{b_0} - u_{b_0}^h\|_{L^\infty(D_{2\epsilon-\eta})} + \|u_{b_1} - u_{b_1}^h\|_{L^\infty(D_{2\epsilon-\eta})} \\
&= \mathcal{O}(\omega^{-5/2}) + \mathcal{O}(\omega^{-3}) + \mathcal{O}(\omega^{-5/2}) \\
&= \mathcal{O}(\omega^{-5/2}).
\end{aligned} \qquad (40)$$

### 4.2. Right-hand side for the far field equation

Given that the Babich's expansion $u_b$ is computed accurately and the cut-off function $\chi_\epsilon$ is known analytically, we can construct the near-field solution $u_{\text{near}} = u_b \chi_\epsilon$. On the other hand, in order to obtain the far-field solution $u_{\text{far}}$, we need to solve the equation (12). Most stability and error analysis for finite element methods [53, 55, 75] rely on the norms of the right-hand side (RHS). Moreover, numerical experiments show that the error of numerical solution is tightly bounded by the norm of RHS in this case, i.e. if the norm of the RHS is $\omega$ dependent, the error will grow as $\omega$ grows. Thus, it is crucial to have asymptotic orders of the norms of the RHS in (12) with respect to $\omega$ and $\epsilon$. We use the analytical expression of the cut-off function and Babich's expansion to obtain such scalings in $L^\infty$ and $L^2$ norm.

| $\Omega \subseteq \mathbb{R}^2$ | $\nabla \chi_\epsilon$ | $\Delta \chi_\epsilon$ | $2\nabla u_b \cdot \nabla \chi_\epsilon$ | $u_b \Delta \chi_\epsilon$ | $2\nabla u_b \cdot \nabla \chi_\epsilon + u_b \Delta \chi_\epsilon$ |
|---|---|---|---|---|---|
| $\|\cdot\|_{L^\infty(\Omega)}$ | $\mathcal{O}(\epsilon^{-1})$ | $\mathcal{O}(\epsilon^{-2})$ | $\mathcal{O}(\omega^{\frac{1}{2}}\epsilon^{-\frac{3}{2}})$ | $\mathcal{O}(\omega^{-\frac{1}{2}}\epsilon^{-\frac{5}{2}})$ | $\mathcal{O}(\omega^{\frac{1}{2}}\epsilon^{-\frac{3}{2}}) + \mathcal{O}(\omega^{-\frac{1}{2}}\epsilon^{-\frac{5}{2}})$ |
| $\|\cdot\|_{L^2(\Omega)}$ | $\mathcal{O}(\epsilon^0)$ | $\mathcal{O}(\epsilon^{-1})$ | $\mathcal{O}(\omega^{\frac{1}{2}}\epsilon^{-\frac{1}{2}})$ | $\mathcal{O}(\omega^{-\frac{1}{2}}\epsilon^{-\frac{3}{2}})$ | $\mathcal{O}(\omega^{\frac{1}{2}}\epsilon^{-\frac{1}{2}}) + \mathcal{O}(\omega^{-\frac{1}{2}}\epsilon^{-\frac{3}{2}})$ |

Table 1: Asymptotic orders of the right-hand side with respect to $\omega$ and $\epsilon$.

We define the smooth cut-off function

$$\chi_\epsilon(\mathbf{x}, \mathbf{x}_0) = \begin{cases} 1, & \text{if } |\mathbf{x} - \mathbf{x}_0| \leq \epsilon, \\ \exp\left(\frac{2e^{-1/t}}{t-1}\right), & \text{if } \epsilon < |\mathbf{x} - \mathbf{x}_0| < 2\epsilon, t = \frac{|\mathbf{x}-\mathbf{x}_0|}{\epsilon} - 1, \\ 0, & \text{if } |\mathbf{x} - \mathbf{x}_0| \geq 2\epsilon, \end{cases} \quad (41)$$

and we choose $\epsilon$ small enough such that $\chi_\epsilon(\mathbf{x}, \mathbf{x}_0)$ is compactly supported within the computational domain $\Omega$, i.e. $D_{2\epsilon} \subseteq \Omega$ and no caustic or multi-pathing has occured. We can easily verify that $\chi_\epsilon \in \mathcal{C}^\infty(\mathbb{R})$, $\nabla \chi_\epsilon$ and $\Delta \chi_\epsilon$ vanishes in $D_{2\epsilon-\epsilon}^c := \{\mathbf{x} \in \Omega : |\mathbf{x} - \mathbf{x}_0| > 2\epsilon \text{ or } |\mathbf{x} - \mathbf{x}_0| < \epsilon\}$. Moreover, based on the derivatives of the first kind Hankel functions (31) and its asymptotic expansions (32), we have the asymptotic orders of the right hand side as displayed in Table 1. The scalings imply that when $\omega \epsilon \gg 1$, the first term $2\nabla u_b \cdot \nabla \chi_\epsilon$ dominates the right-hand side. In the high-frequency regime, we pick small but fixed $\epsilon$. Hence the right-hand side scales as $O(\omega^{\frac{1}{2}})$ as $\omega \to \infty$.

*4.3. Far-field solution: ray-FEM*

After constructing the $\omega$-dependent right hand side of (12), we use ray-FEM to solve the equation with ray directions extracted numerically by NMLA from computed low-frequency wave field. Now we provide an upper bound of the approximation error

$$\inf_{u_h \in V_{Ray}^h(\mathcal{T}_h)} \frac{\|u_{\text{far}} - u_{\text{far}}^h\|_{L^2(\Omega)}}{\|u_{\text{far}}\|_{L^2(\Omega)}}, \quad (42)$$

where the ray-FEM space, $V_{Ray}^h(\mathcal{T}_h)$, incorporates the estimated ray directions $\{\widehat{\mathbf{d}}_j^h\}$. From Appendix A and [34], the error for estimation by NMLA is $|\widehat{\mathbf{d}}_j - \widehat{\mathbf{d}}_j^h| \sim \mathcal{O}(\omega^{-1/2})$.

For the simplicity of error analysis, we assume that there is no ray crossing in the domain $\Omega$, no reflections from the boundary $\partial\Omega$, and the Babich's expansion (13) is the exact solution of equation (1). Then the far-field solution to equation (12) is

$$u_{\text{far}} = u - u_{\text{near}} = u_b - u_b \chi_\epsilon = (1 - \chi_\epsilon) u_b. \tag{43}$$

Note that $(1 - \chi_\epsilon)u_b$ vanishes in the disk $D_\epsilon$, using the first term $u_{b_0}$ to approximate Babich's expansion outside disk $D_\epsilon$ we have truncation error $\mathcal{O}(\omega^{-3/2})$ similar to (34),

$$u_{\text{far}} = (1 - \chi_\epsilon) u_b = (1 - \chi_\epsilon) u_{b_0} + \mathcal{O}(\omega^{-3/2}) = \omega^{-1/2} A(\mathbf{x}) e^{i\omega\phi(\mathbf{x})} + \mathcal{O}(\omega^{-3/2}), \tag{44}$$

where $A(x) = i\frac{\sqrt{\pi\omega}}{2}(1 - \chi_\epsilon) H_0^{(1)}(\omega\phi(\mathbf{x})) e^{i\omega\phi(\mathbf{x})}$ is a smooth amplitude function with a support outside disk $D_\epsilon$. Thus, asymptotically we have

$$\|u_{\text{far}}\|_{L^2(\Omega)} = \mathcal{O}(\omega^{-1/2}). \tag{45}$$

Moreover, we denote by

$$\begin{aligned} u_I &= \sum_{j=1}^{N_h} \omega^{-1/2} A(\mathbf{x}_j) \varphi_j(\mathbf{x}) e^{i\omega[\phi(\mathbf{x}_j) + 1/c(\mathbf{x}_j)\widehat{\mathbf{d}}_j \cdot (\mathbf{x} - \mathbf{x}_j)]}, \\ u_I^h &= \sum_{j=1}^{N_h} \omega^{-1/2} A(\mathbf{x}_j) \varphi_j(\mathbf{x}) e^{i\omega[\phi(\mathbf{x}_j) + 1/c(\mathbf{x}_j)\widehat{\mathbf{d}}_j^h \cdot (\mathbf{x} - \mathbf{x}_j)]} \end{aligned} \tag{46}$$

the nodal interpolations of the solution in $V_{Ray}(\mathcal{T}_h)$ and $V_{Ray}^h(\mathcal{T}_h)$ with exact and numerical ray direction, respectively. For smooth $A$ and $\phi$, from [34] we have

$$\|u_{\text{far}} - u_I\|_{L^2(\Omega)} \lesssim \omega^{-1/2} h^2 |A|_{H^2(\Omega)} + \omega^{1/2} h^2 \|A\|_{L^\infty(\Omega)} \|\nabla^2 \phi\|_{L^\infty(\Omega)} + \mathcal{O}(\omega^{-3/2}), \tag{47}$$

and

$$\|u_I - u_I^h\|_{L^2(\Omega)} \lesssim \omega^{1/2} h \|A\|_{L^\infty(\Omega)} \|c^{-1}\|_{L^\infty(\Omega)} \|\widehat{\mathbf{d}}_j - \widehat{\mathbf{d}}_j^h\|_{L^\infty(\Omega)} \lesssim h \|A\|_{L^\infty(\Omega)} \|c^{-1}\|_{L^\infty(\Omega)}, \tag{48}$$

the constants in $\lesssim$ only depend on the domain $\Omega$. Hence, we have the error estimate,

more compactly with respect to $\omega$ on the mesh with $h = \mathcal{O}(\omega^{-1})$,

$$\begin{aligned}
\inf_{u_{\text{far}}^h \in V_{Ray}^h(\mathcal{T}_h)} \|u_{\text{far}} - u_{\text{far}}^h\|_{L^2(\Omega)} &\leq \|u_{\text{far}} - u_I^h\|_{L^2(\Omega)} \\
&\leq \|u_{\text{far}} - u_I\|_{L^2(\Omega)} + \|u_I - u_I^h\|_{L^2(\Omega)} \\
&= \mathcal{O}(\omega^{-1/2}h^2 + \omega^{1/2}h^2 + h + \omega^{-3/2}) \\
&= \mathcal{O}(\omega^{-1}).
\end{aligned} \tag{49}$$

Therefore,

$$\inf_{u_h \in V_{Ray}^h(\mathcal{T}_h)} \frac{\|u_{\text{far}} - u_{\text{far}}^h\|_{L^2(\Omega)}}{\|u_{\text{far}}\|_{L^2(\Omega)}} = \mathcal{O}(\omega^{-1/2}). \tag{50}$$

We point out that the desirable relative convergence rate in this case is $\mathcal{O}(\omega^{-1})$, which has the same order as the geometric optics ansatz (2). However, the error in the estimation of dominant wave directions using NMLA, which is $\mathcal{O}(\omega^{-1/2})$, dominates the total error and becomes the bottleneck to improve the overall convergence order. This error is due to the deviation of a general wave-front from a plane-wave form which is one of the underpinnings of assumptions for micro local analysis. See Appendix A and [34] for more details.

### 4.4. Adding near-field solution and far-field solution

Adding the near-field solution $u_{\text{near}}$ with the far-field solution $u_{\text{far}}$ and considering the error excluding a small disk $D_\eta$, we can obtain the error estimate for the numerical solution to (1). Indeed,

$$\begin{aligned}
\|u - u_h\|_{L^2(\Omega \setminus D_\eta)} &= \|(u_{\text{near}} + u_{\text{far}}) - (u_{\text{near}}^h + u_{\text{far}}^h)\|_{L^2(\Omega \setminus D_\eta)} \\
&\leq \|(u - u_b^h)\chi_\epsilon\|_{L^2(\Omega \setminus D_\eta)} + \|u_{\text{far}} - u_{\text{far}}^h\|_{L^2(\Omega \setminus D_\eta)} \\
&\leq \|u - u_b^h\|_{L^2(D_{2\epsilon - \eta})} + \|u_{\text{far}} - u_{\text{far}}^h\|_{L^2(\Omega)},
\end{aligned} \tag{51}$$

from (40) and (49), we have

$$\inf_{u_{\text{far}}^h \in V_{Ray}^h(\mathcal{T}_h)} \|u - (u_b^h \chi_\epsilon + u_{smooth}^h)\|_{L^2(\Omega \setminus D_\eta)} = \mathcal{O}(\omega^{-1}). \tag{52}$$

Moreover, based on asymptotic forms in (32) and (33), we have

$$\|u\|_{L^2(\Omega \setminus D_\eta)} = \mathcal{O}(\eta^{-1/2}\omega^{-1/2}), \tag{53}$$

and finally we obtain

$$\inf_{u_{\text{far}}^h \in V_{Ray}^h(\mathcal{T}_h)} \frac{\|u - (u_b^h \chi_\epsilon + u_{\text{far}}^h)\|_{L^2(\Omega \setminus D_\eta)}}{\|u\|_{L^2(\Omega \setminus D_\eta)}} = \mathcal{O}(\eta^{1/2} \omega^{-1/2}), \tag{54}$$

where the constant in $\mathcal{O}(\cdot)$ only depends on constants $C_1$ and $C_2$.

## 5. Fast Linear Solver

The overall complexity claim mentioned in the introduction, is based on the assumption that the linear system provided by both finite element formulations, the standard FEM and the ray-FEM, which we write in a generic form as

$$\mathbf{H}\mathbf{u} = \mathbf{f}, \tag{55}$$

can be solved in linear complexity (up to poly-logarithmic factors).

This can be achieved using a variation of the method of polarized traces [80], which was developed in [34]. The method of itself can be categorized as a domain decomposition method that encompasses the following aspects:

- layered domain decomposition;

- absorbing boundary conditions between subdomains implemented via PML [18];

- transmission conditions issued from a discrete Green's representation formula;

- efficient preconditioner arising from localization of the waves via an incomplete Green's formula.

The method of polarized traces aims to solve the global linear system by solving the local systems, in which the absorbing boundary condition between subdomains, helps to reduce non-physical reflexions due to the truncation of the domain. In order to solve the global system, or in this case to find a good approximate solution, we need to "glue" the subdomains together. The coupling between subdomains is achieved by using the discrete Green's representation formula as explained in [34]. The coupling naturally leads to a smaller integral system posed on the interfaces between subdomains, which can be easily preconditioned by localizing the waves using the incomplete Green's formula. As in [34] we reduce the off-line cost by using a matrix-free formulation (see Chapter 2 in [78]).

In general, the number of iterations for convergence will depend on the quality of the absorbing boundary conditions, and the wave speed. In the best case, the number of iterations will depend on the number of physical reflections across subdomains. For a smooth and fixed wave speed, several numerical experiments indicate that the number of iterations to convergence is weakly dependent on the frequency; i.e., the number of iterations scales as $\mathcal{O}(\log \omega)$, meaning that the cost is dominated by the factorization and solve of each local linear system.

If we suppose that $\Omega$ is discretized into $N = \mathcal{O}(n^2)$ elements, and that each slab is only $\mathcal{O}(1)$ elements thick, then we have that factorizing all the local problems using a multifrontal method [26, 37] (in this case UMFPACK [24]) will have an asymptotic cost of $\mathcal{O}(n)$, which has to be performed $\mathcal{O}(n)$ times, leading to an off-line cost of $\mathcal{O}(N)$.

For the preconditioning, we need to solve $\mathcal{O}(n)$ quasi 1-D linear systems, which can be performed in $\mathcal{O}(n)$ time, leading to a linear complexity for each iteration. This, however, depends on the eventual growth of the auxiliary degrees of freedom corresponding to the PMLs. As it will be shown below, for the low-frequency problem, we need to increase the number of PML points as $\mathcal{O}(\log \omega)$, to maintain the same convergence rate, which is normally achieved in $\mathcal{O}(\log \omega)$ iterations. Thus, the overall complexity is linear up to poly-log factors.

## 6. Algorithms

In this section we provide Algorithm 3 to implement the hybrid method for solving the high-frequency problem (1) with point source terms. The full algorithm starts by computing low-frequency wave field and ray directions in the near field using Babich's expansion in line 3. Then it solves a relative low-frequency Helmholtz equation (30) with the standard FEM to probe the medium and to learn the dominant ray directions in the far field within line 5-8. Afterwards, the high-frequency wave field near the source is computed using the Babich's expansion in line 9. Moreover, the far-field solution of the high-frequency equation (12) is computed using the ray-FEM which incorporates those learned local dominant ray directions in line 11. Lastly, the ray estimation is improved by iteratively applying NMLA to high-frequency waves in the far field and then the resulting rays are incorporated into the ray-FEM space to obtain more accurate high-frequency waves in the while loop from line 14 to line 21.

**Remark 3.** *We refer the reader to [34] for details about* S-FEM*, RayLearning and* Ray-FEM *algorithms. Here $g_{\boldsymbol{u}_{\omega}^{far}}$ and $g_{u_{\omega}^{far}}$ are PML boundary data. In our numerical experiments, we only need one or two iterations in the while loop.*

**Algorithm 3** Hybrid High-Frequency Helmholtz Solver

1: **function** $\mathbf{u}_{\omega,h} = \text{Hybrid-Solver}(\mathbf{x}_0, c, \omega)$
2:     $\widetilde{\omega} \sim \sqrt{\omega}$, $h \sim \omega^{-1}$, $h_c \sim \omega^{-\frac{1}{2}}$     ▷ Low-frequency and mesh sizes
3:     $[u_{\widetilde{\omega},b}^h, \nabla u_{\widetilde{\omega},b}^h, \{\boldsymbol{d}_b^h\}_{D_{2\epsilon}}] = Babich(\mathbf{x}_0, c, \widetilde{\omega}, h)$   ▷ Babich's expansion in disk $D_{2\epsilon}$
4:     $f_{\widetilde{\omega}} = 2\nabla u_{\widetilde{\omega},b}^h \cdot \nabla \chi_\epsilon + u_{\widetilde{\omega},b}^h \Delta \chi_\epsilon$     ▷ Right hand side in (12)
5:     $\mathbf{u}_{\widetilde{\omega},h}^{far} = S\text{-}FEM(\widetilde{\omega}, h, c, f_{\widetilde{\omega}}, g_{\mathbf{u}_{\widetilde{\omega}}^{far}})$     ▷ Low-frequency solution to (30)
6:     $\mathbf{u}_{\widetilde{\omega},h} = \mathbf{u}_{\widetilde{\omega},h}^{far} + u_{\widetilde{\omega},b}^h \chi_\epsilon$     ▷ Low-frequency waves
7:     $\{\boldsymbol{d}_{\widetilde{\omega}}^h\}_{\Omega \setminus D_{2\epsilon}} = RayLearning(\widetilde{\omega}, h, h_c, c, \mathbf{u}_{\widetilde{\omega},h})$   ▷ Ray learning in far field
8:     $\{\boldsymbol{d}_{\widetilde{\omega}}^h\}_\Omega = \{\boldsymbol{d}_{\widetilde{\omega}}^h\}_{\Omega \setminus D_{2\epsilon}} \cup \{\boldsymbol{d}_b^h\}_{D_{2\epsilon}}$     ▷ Low-frequency rays
9:     $[u_{\omega,b}^h, \nabla u_{\omega,b}^h] = Babich(\mathbf{x}_0, c, \omega, h)$
10:    $f_\omega = 2\nabla u_{\omega,b}^h \cdot \nabla \chi_\epsilon + u_{\omega,b}^h \Delta \chi_\epsilon$
11:    $\mathbf{u}_{\mathbf{d}_{\widetilde{\omega}},h}^{far} = Ray\text{-}FEM(\omega, h, c, f_\omega, g_{u_\omega^{far}}, \{\boldsymbol{d}_{\widetilde{\omega}}^h\}_\Omega)$  ▷ High-frequency solution to (12)
12:    $\mathbf{u}_{\omega,h}^1 = \mathbf{u}_{\mathbf{d}_{\widetilde{\omega}},h}^{far} + u_{\omega,b}^h \chi_\epsilon$     ▷ High-frequency waves
13:    $tol = 1$, $niter = 0$,
14:    **while** $tol > \epsilon$ and $niter < max\_iter$ **do**
15:       $\{\mathbf{d}_\omega^h\}_{\Omega \setminus D_{2\epsilon}} = RayLearning(\omega, h, h_c, c, \mathbf{u}_{\omega,h}^1)$
16:       $\{\boldsymbol{d}_\omega^h\}_\Omega = \{\boldsymbol{d}_\omega^h\}_{\Omega \setminus D_{2\epsilon}} \cup \{\boldsymbol{d}_b^h\}_{D_{2\epsilon}}$     ▷ High-frequency rays
17:       $\mathbf{u}_{\mathbf{d}_\omega,h}^{far} = Ray\text{-}FEM(\omega, h, c, f_\omega, g_{u_\omega^{far}}, \{\mathbf{d}_\omega^h\}_\Omega)$
18:       $\mathbf{u}_{\omega,h}^2 = \mathbf{u}_{\mathbf{d}_\omega,h}^{far} + u_{\omega,b}^h \chi_\epsilon$
19:       $tol = \|\mathbf{u}_{\omega,h}^1 - \mathbf{u}_{\omega,h}^2\|_{L^2(\Omega)} / \|\mathbf{u}_{\omega,h}^2\|_{L^2(\Omega)}$
20:       $niter = niter + 1$, $\mathbf{u}_{\omega,h}^1 = \mathbf{u}_{\omega,h}^2$
21:    **end while**
22:    $\mathbf{u}_{\omega,h} = \mathbf{u}_{\omega,h}^1$
23: **end function**

We provide a succinct complexity analysis of the full algorithm 3 in terms of $\omega$, which is summarized in Table 2. The overall complexity includes:

- the complexity to compute the Babich's expansion,
- the complexity of learning ray directions by NMLA, and
- the complexity of the linear solver for the discretized systems from both the standard FEMs (in the low-frequency case), and the ray-FEMs (for high-frequency problem) Helmholtz equations.

It is well studied in Section 5 and [34] that the ray learning stage and fast solvers for the linear systems for both low-frequency and high-frequency problems have linear complexity $\mathcal{O}(\omega^d)$ up to some poly-log factors, as depicted in Table 2. The only extra cost we need to analyze is the cost of computing the Babich's expansion in Algorithm 1 in the uniformly discretized mesh with $h = \mathcal{O}(\omega^{-1})$, which implies $N = \mathcal{O}(\omega^d)$ grid points in the computational domain.

According to [59, 46], by using the high-order LxF-WENO schemes to compute the Babich's ingredients in Algorithm 2, the computational complexity is $\mathcal{O}(N \log N)$, when those asymptotic ingredients are applied in Algorithm 1 to construct the Babich's expansion with linear complexity $\mathcal{O}(N)$. Consequently, the overall complexity for computing Babich's expansion is $\mathcal{O}(\omega^d \log \omega)$ and the overall complexity for the whole hybrid solver is $\mathcal{O}(\omega^d)$ up to a poly-log factors.

| Methods | Babich | S-FEM | Learning | Ray-FEM | Hybrid solver |
|---|---|---|---|---|---|
| Frequency | $\sqrt{\omega}$ or $\omega$ | $\sqrt{\omega}$ | $\sqrt{\omega}$ or $\omega$ | $\omega$ | $\omega$ |
| Complexity | $\mathcal{O}(\omega^d \log \omega)$ | $\mathcal{O}(\omega^d \log^3 \omega)$ | $\mathcal{O}(\omega^d)$ | $\mathcal{O}(\omega^d \log \omega)$ | $\mathcal{O}(\omega^d \log^3 \omega)$ |

Table 2: Overall computational complexities with respect to $\omega$ given that the mesh size scaled as $h = \mathcal{O}(\omega^{-1})$.

## 7. Numerical experiments

In this section we provide several numerical examples to test the proposed method and validate our claims. For all cases, the domain of interest $\Omega$ is discretized using a standard triangular mesh with absorbing boundary conditions implemented via PML while varying the wave speed profile and the source term. The first three cases have a unit square domain $\Omega = (-0.5, 0.5)^2$ and the fourth case has a rectangular domain $\Omega = (-1.5, 1.5) \times (-0.5, 0.5)$. The mesh size $h$ is chosen such that the number of grid points per wavelength (NPW) is fixed, i.e., $\omega h = \mathcal{O}(1)$. Moreover, we fix $\epsilon = \frac{1}{2\pi}$ for frequencies $\omega \geq 100\pi$ so that the $L^2$ norm of the right-hand side of (12) is scaled asymptotically as $\mathcal{O}(\omega^{\frac{1}{2}})$ in this frequency regime (see the explanation at the end of Section 4.2).

We use a high-order Gaussian quadrature rule [27] to compute the integrals required to assemble the mass and stiffness matrices in (25), the right hand side in (26), and the relative $L^2$ errors of the ray-FEM solutions. The algorithms are implemented in MATLAB 2015b using the iFEM package [22] to generate the mesh and the finite element matrices. The numerical experiments were executed in a dual socket server with 2 Intel Xeon E5-2670 and 384 GB of RAM.

### 7.1. Homogeneous medium with exact and numerical rays

We compute the numerical solution to the Helmholtz equation (1) in a homogeneous medium, $c(\mathbf{x}) \equiv 1$, with the exact solution given by

$$u_{ex}(x, y) = \frac{i}{4} H_0^{(1)}(\omega \sqrt{x^2 + y^2}). \tag{56}$$

#### 7.1.1. Convergence

Since the Babich's expansion in a homogeneous medium is exactly the first Hankel function, we use the analytical $u_b$ and $\nabla u_b$ to construct the right-hand side of (12), and we check the convergence rate for the far-field solution $u_{\text{far}}$ with both exact and numerically computed (by NMLA) ray directions.

From Section 4.3, if the ray information is known exactly and $h = \mathcal{O}(\omega^{-1})$, then the relative $L^2$ error in the ray-FEM approximation space is $\mathcal{O}(\omega^{-1})$. Fig. 1 left shows that the ray-FEM is stable and it achieves the desired convergence order with fixed NPW.

On the other hand, if the ray information is numerically estimated by NMLA with accuracy order $\mathcal{O}(\omega^{-1/2})$, the optimal approximation error by the ray-FEM is also $\mathcal{O}(\omega^{-1/2})$ [34]. In fact, we use exact radial ray directions in the disk $D_{2\epsilon}$ and numerically learn the ray directions outside this disk. Fig. 1 right indicates the ray-FEM solution with the learned ray information is of the same order $\mathcal{O}(\omega^{-1/2})$.

#### 7.1.2. Complexity

We use the fast solver developed in [34] to solve (12) thus obtaining the far field component of the wave field. From Fig. 2 we can observe that the results are qualitatively equal to the ones obtained in [34], the complexity is linear up to polylog factors. We point out that the complexity is higher for the low-frequency case given that we need to increase the number of PML points as $\mathcal{O}(\log^2 \omega)$ in order to obtain a very mild growth in the number of iterations. We remark that for this case the largest number of waves we have computed is around 500 in each direction.

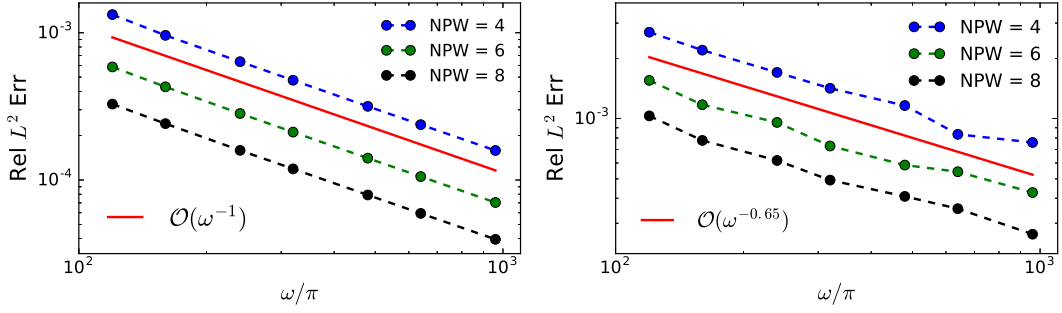

Figure 1: Relative $L^2$ error of smooth part solution to equation (12) for one point-source problem in homogeneous medium, NPW is fixed. Left: exact rays. Right: numerical rays estimated by NMLA.

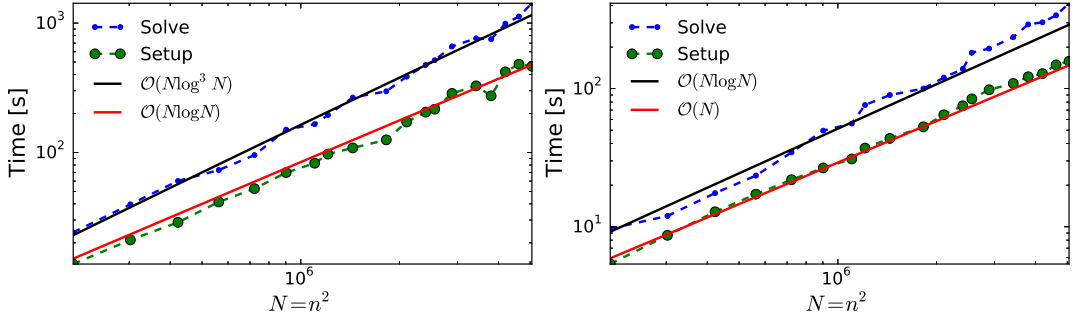

Figure 2: Runtime for solving the Helmholtz equation with a homogeneous wave speed using GMRES preconditioned with the method of polarized traces. The tolerance was set up to $10^{-9}$. Left: runtime for solving the low-frequency problem. Right: Runtime for solving the high-frequency problem with the adaptive basis.

## 7.2. Lippmann-Schwinger equation

We test our algorithm when the wave speed is constant up to a compactly supported heterogeneity. In this case we compute the reference solution by solving the Lippmann-Schwinger equation discretized using the super-algebraically convergent discretization proposed in [1], which is then solved using the fast solver introduced in [81].

In order to use the Lippmann-Schwinger equation, we suppose that the point-source is located far from the heterogeneity. In particular, we set the point-source to be located at $\mathbf{x}_0 = (-0.2, -0.2)$ and the squared slowness to be

$$m(\mathbf{x}) = 1 + 0.2 h(\mathbf{x}, \alpha, \beta) \exp\left(-\frac{r^2(\mathbf{x})}{2\sigma^2}\right), \tag{57}$$

where $\alpha = 0.16, \beta = 0.22, \sigma = 0.15$, $\mathbf{x}_1 = (0.2, 0.2)$, $r(\mathbf{x}) = |\mathbf{x} - \mathbf{x}_1|$, $t(\mathbf{x}, \alpha, \beta) = \frac{r(\mathbf{x}) - \alpha}{\beta - \alpha}$, $P(t) = \frac{2e^{-1/t}}{t-1}$ and

$$h(\mathbf{x}) = \begin{cases} 1, & \text{if } r(\mathbf{x}) \leq \alpha, \\ \exp\left(P(t(\mathbf{x}, \alpha, \beta))\right), & \text{if } \alpha < r(\mathbf{x}) < \beta, \\ 0, & \text{if } r(\mathbf{x}) \geq \beta. \end{cases}$$

In this case, in the disk $D_{2\epsilon}$, the medium is homogeneous so that we can build the right-hand side and the rays analytically; outside of this disk, the medium is heterogeneous and we apply NMLA to estimate the ray directions. Besides, we use Algorithm 3 to compute the ray-FEM solution to the far field equation (12) and compare it to the reference solution described in [81]. Fig. 3 shows that the relative error in the $L^2(\Omega)$ norm follows the desired convergence rate $\mathcal{O}(\omega^{-1/2})$ with fixed NPW. We mention that for this example the largest number of waves we have computed is around 500 in each direction.

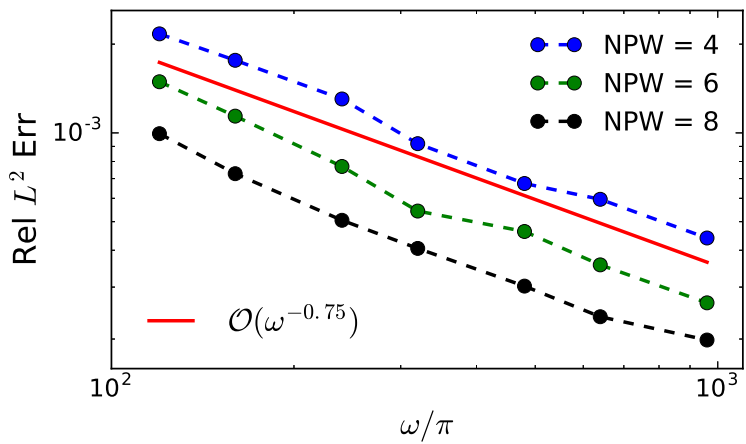

Figure 3: Relative $L^2$ error of smooth part solution to Lippmann-Schwinger equation with squared slowness (57), NPW is fixed.

### 7.3. Wave speed of constant gradient

We provide an example in a heterogeneous medium with wave speed of constant gradient: $c(\mathbf{x}) = c_0 + G_0 \cdot (\mathbf{x} - \mathbf{x}_0)$ with parameters $c_0 = 1$, $G_0 = (0.1, -0.2)$ and $\mathbf{x}_0 = (0, 0)$. The phase function is known analytically [36] and there is no ray crossing in the domain $\Omega = (-0.5, 0.5)^2$. Then Algorithm 1 can produce an accurate solution to (1) in the whole domain and we treat this solution as the reference solution.

We construct the right-hand side with numerically computed $u_b$ and $\nabla u_b$ in the disk $D_{2\epsilon}$ and then apply Algorithm 3 to compute the numerical solution to (1). We compute the relative $L^2$ error with respect to the reference solution. Fig. 4 left shows that the error scales as $\mathcal{O}(\omega^{-1})$ when we use analytical rays outside the disk $D_{2\epsilon}$; on the other hand, Fig. 4 right shows that the error scales as $\mathcal{O}(\omega^{-1/2})$ when we use numerically computed rays by NMLA instead.

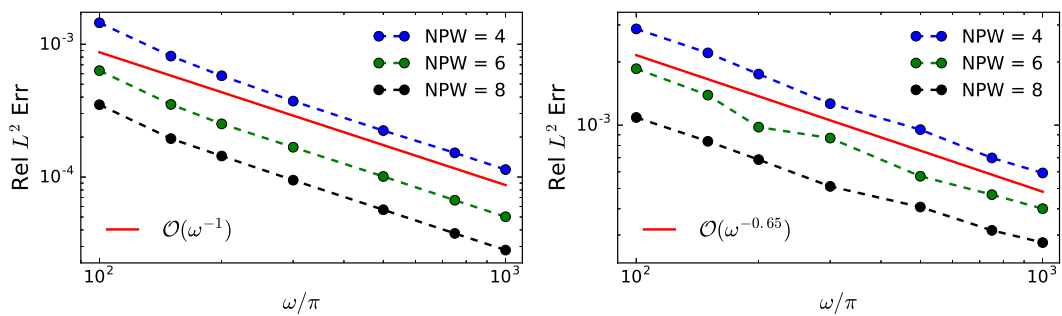

Figure 4: Relative $L^2$ error of numerical solution to Helmholtz equation (1) with constant gradient of velocity, NPW is fixed. Left: analytical rays. Right: numerical rays estimated by NMLA.

### 7.4. Marmousi model

Finally, we apply our method to the Marmousi2 model [52]. Fig. 5 shows the wavespeed, which is smoothed by a convolution with a Gaussian kernel with standard deviation of 100 meters. In this model, we re-scale the computational domain to $\Omega = (-1.5, 1.5) \times (-0.5, 0.5)$ and we locate the point source such that the wave speed is constant in the neighborhood $D_{2\epsilon}$. Within this neighborhood, caustics do not occur so that the Babich's expansion is reduced to the Hankel function, which can compute the wave field and ray directions very accurately. However, in the far field $\Omega \backslash D_{2\epsilon}$, where ray crossing happens and caustics occur, we utilize NMLA to capture only the local dominant ray directions. We select at most four significant ray directions by sorting amplitudes. In addition, we select rays that are well separated with an angle difference at least 15 degrees.

Then we use Algorithm 3 to compute the wave field at 18.75 [Hz] on the mesh with 4 grid points for the smallest wavelength (NPW = 4). The real part of the wave field is shown in Fig. 6. Furthermore, we regard the solution on the mesh with NPW = 16 as the reference solution $u_{ref}$ and compute numerical solutions $u_h$ on different coarser meshes to show the $h$ convergence rate in Table 3. A higher frequency case at 75 [Hz] is shown in Fig. 7. At the highest frequency the solution has roughly 250 wavelengths in the vertical direction and 750 in the horizontal direction.

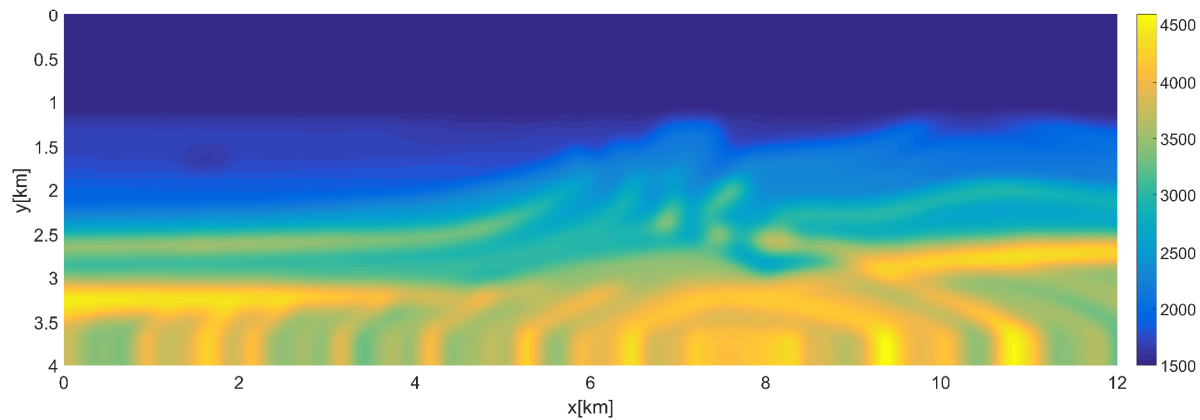

Figure 5: Smoothed Marmousi wave speed model.

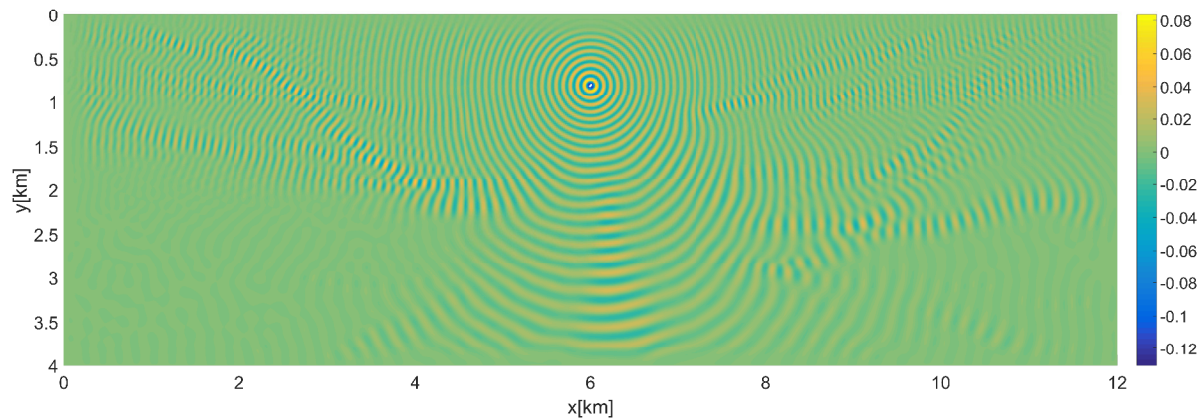

Figure 6: Real part of wave field generated by a point-source at 18.75 [Hz] with NPW = 4 for the smoothed Marmousi model.

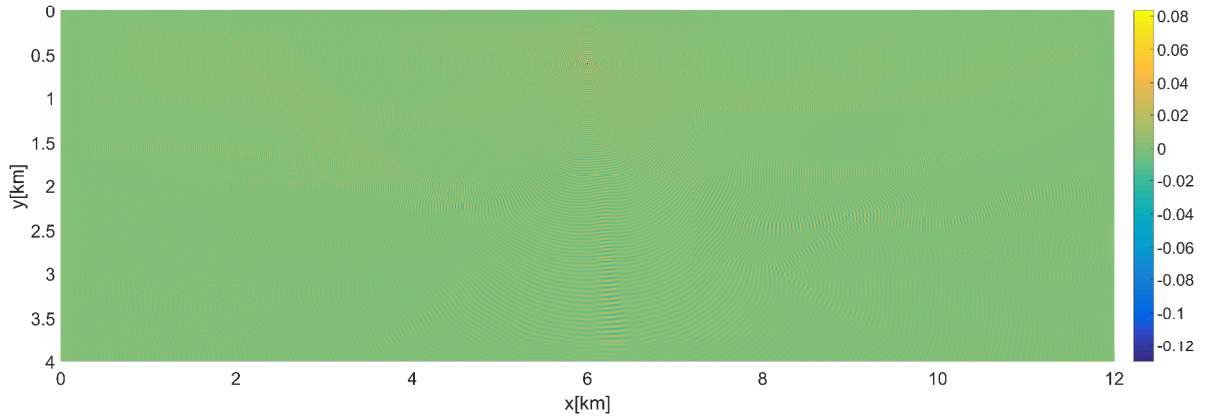

Figure 7: Real part of wave field generated by a point-source at 75 [Hz] with NPW = 4 for the smoothed Marmousi model.

| NPW | 1 | 2 | 4 | 8 |
|---|---|---|---|---|
| $\frac{\|u_h - u_{ref}\|_{L^2(\Omega)}}{\|u_{ref}\|_{L^2(\Omega)}}$ | 0.9645 | 0.2846 | 0.0806 | 0.0218 |
| Order | | 1.76 | 1.82 | 1.89 |

Table 3: Marmousi model $h$ convergence rate.

## 8. Conclusion

We present a hybrid method, which combines the asymptotic Babich's expansion and the ray-FEM method, for solving the high-frequency Helmholtz equation in smooth media with point source terms. The method removes the singularity efficiently by using the Babich's asymptotic expansion in the near field, and computes the far-field solution by the ray-FEM that incorporates local dominant plane-wave directions adaptive to the medium and source. Numerical tests suggest that the method only requires a fixed number of points per wavelength with an asymptotic convergence $\mathcal{O}(\omega^{-1/2})$. The proposed method is coupled with a fast sweeping-type solver achieving an empirical complexity $\mathcal{O}(\omega^2)$ up to a ploy-log factors in 2D.

## Acknowledgement

Zhao is partially supported by NSF Grants (1418422 and 1622490). Qian is partially supported by NSF Grants (1522249 and 1614566).

## Appendix A. NMLA

Here is a brief introduction to NMLA developed in [16, 17]. We suppose that a wave field is a locally weighted superposition of plane-waves having the same wave number and propagating in different directions. The aim of NMLA is to extract the directions and the weights by sampling and processing the wave field locally.

Suppose that a wave field, denoted by $u(\mathbf{x})$, is composed of $N$ plane-waves around an observation point $\mathbf{x}_0$,

$$u(\mathbf{x}) = \sum_{n=1}^{N} B_n e^{ik(\mathbf{x}-\mathbf{x}_0)\cdot \widehat{\mathbf{d}}_n}, \quad |\widehat{\mathbf{d}}_n| = 1. \tag{A.1}$$

We suppose that we can sample the wave field, $u(\mathbf{x})$, and its derivative on a circle $S_r(\mathbf{x}_0)$ centered at $\mathbf{x}_0$ with radius $r$. The wave field can be written under the model assumption in (A.1) as

$$u(\mathbf{x}_0 + r\widehat{\mathbf{s}}) = \sum_{n=1}^{N} B_n e^{i\alpha \widehat{\mathbf{s}}\cdot \widehat{\mathbf{d}}_n}, \quad \alpha = kr, \ \widehat{\mathbf{s}} \in \mathbb{S}^1. \tag{A.2}$$

Furthermore, we define the angle variables $\theta = \theta(\widehat{\mathbf{s}})$ and $\theta_n = \theta(\widehat{\mathbf{d}}_n)$ such that $\widehat{\mathbf{s}} = (\cos\theta, \sin\theta)$, $\widehat{\mathbf{d}}_n = (\cos\theta_n, \sin\theta_n)$, and $\mathbf{x}(\theta) = \mathbf{x}_0 + r\widehat{\mathbf{s}}(\theta)$. Using the angle based notation we sample the impedance quantity on the circle $S_r(\mathbf{x}_0)$,

$$U(\theta) := \frac{1}{ik}\partial_r u(\mathbf{x}(\theta)) + u(\mathbf{x}(\theta)), \tag{A.3}$$

which removes any possible ambiguity due to resonance [16] and improves the robustness to noise for solutions to the Helmholtz equation. We apply the filtering operator $\mathcal{B}$ to $U(\theta)$ defined as

$$\mathcal{B}U(\theta) := \frac{1}{2L_\alpha + 1} \sum_{l=-L_\alpha}^{L_\alpha} \frac{(\mathcal{F}U)_l e^{il\theta}}{(-i)^l(J_l(\alpha) - iJ'_l(\alpha))}, \tag{A.4}$$

where $L_\alpha = \max(1, [\alpha], [\alpha + (\alpha)^{\frac{1}{3}} - 2.5])$, $J_l$ is the Bessel function of order $l$, $J_l'$ is its derivative and

$$(\mathcal{F}U)_l := \frac{1}{2\pi} \int_0^{2\pi} U(\theta) e^{-il\theta} d\theta \tag{A.5}$$

is the $l$-th Fourier coefficient of $U$. It is shown in [16] that

$$\mathcal{B}U(\theta) = \sum_{n=1}^{N} B_n S_{L_\alpha}(\theta - \theta_n), \tag{A.6}$$

where $S_L(\theta) = \frac{\sin([2L+1]\theta/2)}{[2L+1]\sin(\theta/2)}$. As a consequence, we have that if $\alpha = kr \to \infty$ then

$$\lim_{\alpha \to \infty} \mathcal{B}U(\theta) = \begin{cases} B_n, & \text{if } \theta = \theta_n \text{ (or } \widehat{\mathbf{s}} = \widehat{\mathbf{d}}_n \text{ )}; \\ 0, & \text{otherwise}. \end{cases} \tag{A.7}$$

We obtain the directions by picking the peaks in the filtered data in (A.6), and the amplitudes by solving a least square problem with the direction obtained.

However, for applications, the measured data are never a perfect superposition of plane-waves introducing errors in the estimation. We summarize the stability result and error estimate with noisy data from [16, 34] in the sequel.

For simplicity we use the single wave case, i.e., $N = 1$. Moreover, we assume that the measured datum is a perturbation to the perfect plane-wave datum of the form $U(\theta) = U^{plane}(\theta) + \delta U(\theta)$, where $U^{plane}$ denotes a single plane-wave datum in the form of (A.2). Let $\theta^*$ denote the angle for which $\theta \mapsto \mathcal{B}U(\theta)$ is maximum. Assuming that the noise level satisfies

$$||\delta U||_{L^\infty} < \frac{1}{4B^*} |B_1|, \tag{A.8}$$

where $B^* \leq 1$ is a pure constant independent of $\omega$ and $B_1$ is the complex amplitude of the plane-wave. Then the error in the angle estimation is given by

$$|\theta_1 - \theta^*| \leq \frac{2\pi}{2L_\alpha + 1} \sim \mathcal{O}(\frac{1}{\alpha}), \quad \alpha = kr \sim \infty. \tag{A.9}$$

Similar results can be derived for multiple waves $N > 1$. We remark that $\frac{1}{4B^*} \geq 0.25$, which implies that if the relative noise level does not surpass 25% the direction will be estimated within an error of order $\mathcal{O}(\frac{1}{kr})$. In our problem, when NMLA is used to estimate local ray directions, the perturbation is due to the fact that a general wave field is a superposition of curved wave-fronts. In [34], it was shown that the optimal choice of the radius of the observation circle is $r \sim \omega^{-1/2}$ for wave-fronts with bounded curvature, and one can achieve the following error estimate

$$|\theta_1 - \theta^*| = \mathcal{O}(\omega^{-1/2}). \tag{A.10}$$